\documentclass{article}
\usepackage{amsmath}

\title{SOME RATIONAL SOLUTIONS TO PAINLEVE' VI}
\author{Gert Almkvist}
\begin{document}
\maketitle

%\Setdate{}
%\TitlePage{}

Much has been written about rational solutions to Painleve' VI (denoted by $%
P_{VI}(\alpha ,\beta ,\gamma ,\delta )$ )
\[
\]
\[
y^{\prime \prime }=\frac 12(\frac 1y+\frac 1{y-1}+\frac 1{y-x})\text{ }%
y^{\prime }\text{ }^2-(\frac 1x+\frac 1{x-1}+\frac 1{y-x})\text{ }y^{\prime
}+ 
\]
\[
\frac{y(y-1)(y-x)}{x^2(x-1)^2}\left\{ \alpha +\frac{\beta x}{y^2}+\frac{%
\gamma (x-1)}{(y-1)^2}+\frac{\delta x(x-1)}{(y-x)^2}\right\} 
\]
\[
\]
See [1],[2],[3],[4],[5],[6],[7]. Here we give a selfcontained elementary
treatment of how to find two 2-parametric families of rational solutions.In
particular we establish a direct connection between $P_{VI}$ and Jacobi
polynomials (studied in [4]). If Theorem 4.2 in Yuan-Li [7] were true , then
these solutions would be the only rational solutions. But the proof of
Theorem 4.2 in [7] is not correct ( e.g. $R^2$ in formula $(4.4)$ should be $%
\nu ^2R^2$ ). We start with the following 
\[
\]

\textbf{Lemma 1:} Assume that $y$ satisfies the Riccati equation 
\[
\]
\[
x(x-1)\text{ }y^{\prime }=ay^2+(bx+c)y+dx 
\]
\[
\]
where 
\[
\]
\[
a+b+c+d=0 
\]
Then $y$ satisfies 
\[
\]
\[
P_{VI}(\frac{a^2}2,-\frac{d^2}2,\frac{(b+d)^2}2,\frac{1-(c+d+1)^2}2) 
\]
\[
\]

\textbf{Proof:} Differentiate the Riccati equation 
\[
\]
\[
x(x-1)y^{\prime \prime }+(2x-1)y^{\prime }=2ayy^{\prime }+(bx+c)y^{\prime
}+by+d 
\]
\[
\]
Solve for $y^{\prime \prime }$%
\[
\]
\[
y^{\prime \prime }=\frac{(2ay+(b-2)x+c+1)(ay^2+(bx+c)y+dx)}{x^2(x-1)^2}+%
\frac{by+d}{x(x-1)} 
\]
\[
\]
If we substitute 
\[
\]
\[
y^{\prime }=\frac{ay^2+(bx+c)y+dx}{x(x-1)} 
\]
\[
\]
and 
\[
\]
\[
d=-a-b-c 
\]
\[
\]
in 
\[
\frac 12(\frac 1y+\frac 1{y-1}+\frac 1{y-x})\text{ }y^{\prime }\text{ }%
^2-(\frac 1x+\frac 1{x-1}+\frac 1{y-x})\text{ }y^{\prime }+ 
\]
\[
\]
\[
\frac{y(y-1)(y-x)}{x^2(x-1)^2}\left\{ \frac{a^2}2-\frac{d^2x}{2y^2}+\frac{%
(b+d)^2(x-1)}{2(y-1)^2}+\frac{(1-(c+d+1)^2)x(x-1)}{2(y-x)^2}\right\} 
\]
\[
\]
then (using Maple) we see that the right hand sides agree. 
\[
\]
\[
\]

Now we linearize 
\[
\]
\[
x(x-1)\text{ }y^{\prime }=ay^2+(bx+c)y+dx 
\]
\[
\]
by the substitution 
\[
\]
\[
y=-\frac{x(x-1)w^{\prime }}{aw} 
\]
\[
\]
We get 
\[
\]
\[
x(x-1)^2y^{\prime \prime }-(x-1)(rx+s)\text{ }w^{\prime }+tw=0 
\]
\[
\]
where 
\[
\]
\[
r=b-2 
\]
\[
s=c+1 
\]
\[
t=ad 
\]
\[
\]
This is a degenerated Heun equation. It has many polynomial solutions. 
\[
\]

\textbf{Example:}Let 
\[
\]
\[
r=8 
\]
\[
s=2 
\]
\[
t=30 
\]
\[
\]
Then 
\[
\]
\[
a=-6 
\]
\[
b=10 
\]
\[
c=1 
\]
\[
d=-5 
\]
\[
\]
We obtain the solution 
\[
\]
\[
w=1-15x+90x^2-295x^3+594x^4-771x^5+650x^6-345x^7+105x^8-14x^9 
\]
\[
\]
But 
\[
\]
\[
y=\frac{x(x-1)w^{\prime }}{6w} 
\]
\[
\]
collapses to 
\[
\]
\[
y=\frac 12\frac{x(42x^3-70x^2+35x-5)}{(2x-1)(7x^2-7x+1)} 
\]
\[
\]
which solves 
\[
\]
\[
P_{VI}(18,-\frac{25}2,\frac{25}2,-4) 
\]
\[
\]
The collaps depends on the factorization 
\[
\]
\[
w=(1-x)^6(1-9x+21x^2-14x^3) 
\]
\[
\]

Therefore we make the substitution 
\[
\]
\[
w=(1-x)^k\text{ }u 
\]
\[
\]
It follows 
\[
\]
\[
w^{\prime }=(1-x)^k\text{ }u^{\prime }-k(1-x)^{k-1}u 
\]
\[
\]
\[
w^{\prime \prime }=(1-x)^ku^{\prime \prime }-2k(1-x)^{k-1}u^{\prime
}+k(k-1)(1-x)^{k-2}u 
\]
\[
\]
This gives after canceling $(1-x)^{k-2}$%
\[
\]
\[
x(1-x)^2u^{\prime \prime }+(1-x)(rx-2kx+s)u^{\prime }+\left\{
k(k-1)x+t-k(rx+s)\right\} u=0 
\]
\[
\]
So far we still have the factor $x(1-x)^2$ in front of $u^{\prime \prime }$
but by choosing $k$ we can make the coefficient of $u$ divisible by $1-x.$
Hence 
\[
\]
\[
(k(k-1)-kr)x+t-ks= 
\]
\[
\]
\[
(k^2-k-kr)(x-1)+k^2-k-rk+t-sk 
\]
\[
\]
Now choose $k$ such that 
\[
\]
\[
k^2-(r+s+1)k+t=0 
\]
\[
\]
Then we are left with the hypergeometric equation 
\[
\]
\[
x(1-x)u^{\prime \prime }+\left\{ s-(k+(k-r-1)+1)x\right\} u^{\prime
}-k(k-r-1)u=0 
\]
\[
\]
with one solution 
\[
\]
\[
u_1=F(k,k-r-1,s,x) 
\]
\[
\]
where 
\[
\]
\[
F(\alpha ,\beta ,\gamma ,x)=\sum_{n=0}^\infty \frac{(\alpha )_n(\beta )_n}{%
(\gamma )_n}\frac{x^n}{n!} 
\]
\[
\]
Here 
\[
\]
\[
(\alpha )_n=\alpha (\alpha +1)\cdot \cdot \cdot (\alpha +n-1) 
\]
\[
\]
is the Pochhammer symbol. To obtain a rational $y$ we want $u_1$ to be a
polynomial. This is achieved if 
\[
\]
\[
k-r-1=-\mu 
\]
\[
\]
is a negative integer. We express everything in the parameters 
\[
\]
\[
k,\mu ,s 
\]
We get 
\[
\]
\[
r=k+\mu -1 
\]
\[
b=k+\mu +1 
\]
\[
c=s-1 
\]
\[
\]
Since 
\[
\]
\[
a+d=-b-c=-(k+\mu +s) 
\]
\[
ad=t=(r+s+1)\lambda -\lambda ^2=\lambda (\mu +s) 
\]
\[
\]
We get two solutions 
\[
\]
\[
a=-k 
\]
\[
d=-(\mu +s) 
\]
\[
\]
and 
\[
a=-(\mu +s) 
\]
\[
d=-k 
\]
\[
\]

\textbf{First case:} $a=-k$%
\[
\]
Wehave 
\[
\]
\[
y_1=y_1(k,\mu ,s)=\frac{x(x-1)w^{\prime }}{kw}= 
\]
\[
\]
\[
x\left\{ 1+\frac{\mu (1-x)}s\frac{F(k+1,1-\mu ,s+1,x)}{F(k,-\mu ,s)}\right\} 
\]
\[
\]
which solves 
\[
\]
\[
P_{VI}(\frac{k^2}2,-\frac{(\mu +s)^2}2,\frac{(k-s+1)^2}2,\frac{1-\mu ^2}2) 
\]
\[
\]
Here we use that 
\[
\]
\[
\frac d{dx}F(\alpha ,\beta ,\gamma ,x)=\frac{\alpha \beta }\gamma F(\alpha
+1,\beta +1,\gamma +1,x) 
\]
\[
\]

\textbf{Second case: }$a=-(\mu +s).$%
\[
\]
We get 
\[
\]
\[
y_2=y_2(k,\mu ,s)=\frac{x(x-1)w^{\prime }}{(\mu +s)w}= 
\]
\[
\]
\[
\frac{kx}{(\mu +s)}\left\{ 1+\frac{\mu (1-x)}s\frac{F(k+1,1-\mu ,s+1,x)}{%
F(k,-\mu ,s,x)}\right\} 
\]
\[
\]
which solves 
\[
\]
\[
P_{VI}(\frac{(\mu +s)^2}2,-\frac{k^2}2,\frac{(\mu +1)^2}2,\frac{1-(k-s)^2}2) 
\]
\[
\]

These solutions $y_1$ and $y_2$ are rational if $\mu $ is a positive
integer. But you can express $y_2$ in $y_1$ with different parameters since 
\[
\]
\[
y_2(k,\mu ,s)=y_1(\mu +s,k-s,s) 
\]
\[
\]
(but it is not clear that the RHS is rational when $\mu $ is a positive
integer). This equality follows from taking the logatithmic derivative of
the identity 
\[
\]
\[
F(k,-\mu ,s,x)=(1-x)^{\mu +s-k}F(\mu +s,s-k,s,x) 
\]
\[
\]
\[
\]

If $s\notin \mathbf{Z}$ then the hypergeometric equation has a second
solution 
\[
\]
\[
u_2=x^{1-s}F(k-s+1,1-\mu -s,2-s,x)
\]
\[
\]
We note that a nontrivial linear combination of $u_1$ and $u_2$ never can be
a polynomial or have a rational logarithmic derivative. 
\[
\]
\[
\]

\textbf{Case 1:} $a=-k$%
\[
\]
We get 
\[
\]
\[
y_3=y_3(k,\mu ,s)=\frac{x(x-1)}k\frac d{dx}\log \left\{
x^{1-s}(1-x)^kF(k-s+1,1-\mu -s,2-s,x)\right\} = 
\]
\[
\]
\[
\frac{x(x-1)}k\left\{ \frac{1-s}x-\frac k{1-x}+\frac{(k-s+1)(1-\mu -s)}{2-s}%
\frac{F(k-s+2,2-\mu -s,3-s,x)}{F(k-s+1,1-\mu -s,2-s,x)}\right\} 
\]
\[
\]
which solves 
\[
\]
\[
P_{VI}(\frac{k^2}2,-\frac{(\mu +s)^2}2,\frac{(k-s+1)^2}2,\frac{1-\mu ^2}2) 
\]
\[
\]
\[
\]
We see that $y_3$ is rational if one of $1-\mu -s$ or $k-s+1$ is a negative
integer ( and $s\notin \mathbf{Z}$ ). If also $\mu $ is a positive integer
then the above $P_{VI}$ has two rational solutions $y_1$ and $y_3.$%
\[
\]
\[
\]

\textbf{Case 2: }$a=-(\mu +s).$%
\[
\]
We get 
\[
\]
\[
y_4=y_4(k,\mu ,s)=\frac{x(x-1)}{\mu +s}\frac d{dx}\log \left\{
x^{1-s}(1-x)^kF(k-s+1,1-\mu -s,2-s,x)\right\} = 
\]
\[
\]
\[
\frac{x(x-1)}{\mu +s}\left\{ \frac{1-s}x-\frac k{1-x}+\frac{(k-s+1)(1-\mu -s)%
}{2-s}\frac{F(k-s+2,2-\mu -s,3-s,x)}{F(k-s+1,1-\mu -s,2-s,x)}\right\} 
\]
\[
\]
which solves 
\[
\]
\[
P_{VI}(\frac{(\mu +s)^2}2,-\frac{k^2}2,\frac{(\mu +1)^2}2,\frac{1-(k-s)^2}2) 
\]
\[
\]
Note that $y_4$ is rational in the same cases as $y_3.$ Here we also have 
\[
\]
\[
y_4(k,\mu ,s)=y_3(\mu +s,k-s,s) 
\]
\[
\]

\textbf{Remark:} Sometimes you get rational solutions unexpectedly. E.g. we
have 
\[
\]
\[
y_2(4,\sqrt{2},2)=\frac{(3-\sqrt{2})x(7x^2-16x+4x\sqrt{2}+12-6\sqrt{2})}{%
7x^2+6x\sqrt{2}-18x+24-15\sqrt{2}} 
\]
\[
\]
which is equal to 
\[
y_1(2+\sqrt{2},2,2) 
\]
\[
\]
It solves 
\[
P_{VI}(\frac{(2+\sqrt{2})^2}2,-8,\frac{(1+\sqrt{2})^2}2,-\frac 32) 
\]
\[
\]

We collect our results in the following Theorem ( the solutions $y_2$ and $%
y_4$ can be expressed in $y_1$ and $y_3$ respectively) 
\[
\]

\textbf{Theorem:} The equation 
\[
P_{VI}(\frac{k^2}2,-\frac{(\mu +s)^2}2,\frac{(k-s+1)^2}2,\frac{1-\mu ^2}2) 
\]
\[
\]
has the solutions 
\[
\]
(a) 
\[
y=x\left\{ 1+\frac{\mu (1-x)}s\frac{F(k+1,1-\mu ,s+1,x)}{F(k,-\mu ,s,x)}%
\right\} 
\]
\[
\]
It is rational in the following cases 
\[
\]
(i) $\mu $ is a positive integer 
\[
\]
(ii) $k$ is a negative integer 
\[
\]
(iii) $k-s$ is a positive integer 
\[
\]
(iv) $\mu +s$ is a negative integer 
\[
\]
\[
\]
(b) If $s\notin \mathbf{Z}$ there is a second solution 
\[
\]
\[
y=\frac{x(x-1)}k\left\{ \frac{1-s}x-\frac k{1-x}+\frac{(k-s+1)(1-\mu -s)}{2-s%
}\frac{F(k-s+2,2-\mu -s,3-s,x)}{F(k-s+1,1-\mu -s,2-s,x)}\right\} 
\]
\[
\]
\[
\]
It is rational in the following cases 
\[
\]
(i) $1-\mu -s$ is a negative integer 
\[
\]
(ii) $k-s+1$ is a negative integer 
\[
\]
(iii) $k$ is a positive integer 
\[
\]
(iv) $\mu $ is a negative integer 
\[
\]
\[
\]

\textbf{The case }$\mathbf{\alpha =\delta =0.}$%
\[
\]

It follows from Garnier [2] and Gromak-Lukashevich [3] that if 
\[
\]
\[
y^{\prime \text{ }2}=\frac{2(y-x)^2}{x^2(x-1)^2}\left\{ (\beta +\gamma
)y-\beta \right\} 
\]
\[
\]
then $y$ satisfies 
\[
P_{VI}(0,\beta ,\gamma ,0) 
\]
\[
\]
Then make the substitution 
\[
\]
\[
y=\frac \beta {\beta +\gamma }+\frac{2x^2(x-1)^2u^{\prime \text{ }2}}{(\beta
+\gamma )u^2} 
\]
\[
\]
Then we get the linear equation 
\[
\]
\[
u^{\prime \prime }+\frac{2x-1}{x(x-1)}u^{\prime }-\frac 1{x^2(x-1)^2}\left\{ 
\frac{\beta +\gamma }2x-\frac \beta 2\right\} u=0 
\]
\[
\]
To get a hypergeometric equation we put 
\[
\]
\[
u=x^s(x-1)^tv 
\]
\[
\]
Then 
\[
\]
\[
v^{\prime \prime }+\left\{ 2(\frac sx+\frac t{x-1})+\frac{2x-1}{x(x-1)}%
\right\} v^{\prime }+ 
\]
\[
\left\{ \frac{s^2-s}{x^2}+\frac{t^2-t}{(x-1)^2}+\frac{2st}{x(x-1)}+\frac{2x-1%
}{x(x-1)}(\frac sx+\frac t{x-1})-\frac{\beta +\gamma }2x-\frac \beta
2\right\} v=0 
\]
\[
\]
Now we choose $s$ and $t$ such that the numerator of the coefficient of $v$
is divisible by $x(x-1).$ We obtain 
\[
\]
\[
\beta =-2s^2 
\]
\[
\gamma =-2t^2 
\]
\[
\]
and 
\[
x(1-x)v^{\prime \prime }+\left\{ 2s+1-(s+t+(s+t+1)+1)x\right\} v^{\prime
}-(s+t)(s+t+1)v=0 
\]
\[
\]
Hence 
\[
v_1=F(s+t,s+t+1,2s+1,x) 
\]
\[
\]
and 
\[
u_1=x^s(x-1)^tF(s+t,s+t+1,2s+1,x) 
\]
\[
\]
If $2s+1\notin \mathbf{Z}$ then we get a second solution 
\[
\]
\[
u_2=x^{-s}(x-1)^tF(t-s,t-s+1,-2s+1,x) 
\]
\[
\]
so 
\[
u_2(s,t)=u_1(-s,t) 
\]
\[
\]
is nothing new. To simplify notation let 
\[
\]
\[
n=s+t 
\]
\[
r=1+2s 
\]
\[
\]
Then 
\[
u_1=x^{(r-1)/2}(x-1)^{n-(r-1)/2}F(n,n+1,r,x) 
\]
\[
\]

\textbf{Theorem 2: }We have 
\[
\]
\[
y(n,r)=\frac 1{n(n-r+1)}\left\{ 
\begin{array}{c}
-\frac{(r-1)^2}4+ \\ 
\left[ nx-\frac{r-1}2+\frac{n(n+1)x(x-1)}r\frac{F(n+1,n+2,r+1,x)}{%
F(n,n+1,r,x)}\right] ^2
\end{array}
\right\} 
\]
solves 
\[
P_{VI}(0,-\frac{(r-1)^2}2,\frac{(2n-r+1)^2}2,0) 
\]
\[
\]
which is rational if $n$ is a negative integer or if $n-r$ is a nonnegative
integer. 
\[
\]

Now we can give an explicit counterexample to Theorem 4.2 in [7]. Let 
\[
\]
\[
p=0 
\]
\[
\lambda =6 
\]
\[
r=10 
\]
\[
q=1 
\]
\[
C=-2 
\]
\[
\]
in Theorem D in [7]. Then 
\[
w=\frac x4\frac{(x+8)(x^2+14x+21)}{(2x+7)^2} 
\]
\[
\]
solves 
\[
P_{VI}(0,-18,50,0) 
\]
\[
\]
but 
\[
v=\frac{4(2x+7)^2}{(x+7)(x+8)(x^2+7x+28)} 
\]
\[
\]
is not zero.

\subsection*{References:} 
\begin{enumerate}
\item A.S.Fokas, M.J.Ablowitz, On a unified approach to transformations and
elementary solutions to Painlev\'{e} equations, J. Math. Phys. 23 (1982),
2033-2043. 

\item R.Garnier, Contribution a l'etude des solutions de l'equation (V) de
Painlev\'{e}. 

\item V.I.Gromak, N.A.Lukashevich, Special classes of solutions of Painlev\'{e}
equations, Diff Eqns. 18 (1982), 317-326. 

\item L.Haine, J-P.Semengue, The Jacobi polynomial ensemble and the
Painlev\'{e}VI equation, J. Math. Phys. 40 (1999), 2117-2134. 

\item M.Mazzocco, Rational solutions of the Painlev\'{e} VI equation,
arXiv:nlin SI/0007036. 

\item K.Okamoto, Studies on the Painlev\'{e} equations I. Sixth Painlev\'{e}
equation, Ann. Mat. Pura Appl. 146 (1987), 337-381. 

\item Yuan Wenjun, Li Yezhou, Rational solutions of Painlev\'{e} equations,
Canad. J. Math. 54 (2002), 648-670. 
\end{enumerate}

\vfill

University of Lund

Box 118

S22100 Lund

Sweden

gert@maths.lth.se

\end{document}